\documentclass[14pt,a4paper]{article}

\usepackage[utf8]{inputenc}
\usepackage[T2A]{fontenc}
\usepackage[english,russian]{babel}

\usepackage{indentfirst} 
\usepackage[14pt]{extsizes}
\usepackage{amssymb}
\usepackage{amsmath}
\usepackage{amsthm}
\usepackage{pifont}
\usepackage{graphicx}
\usepackage[colorinlistoftodos]{todonotes}
\usepackage[colorlinks=true, allcolors=blue]{hyperref}
\usepackage{algorithm}
\usepackage{algorithmic}

\newtheorem{theorem}{Теорема}

\def\UDK#1{{\leftline{УДК {#1}}}}

\def\aig#1{{\color{black}#1}}
\def\ag#1{{\color{black}#1}}

\usepackage{geometry}
\begin{document}
\renewcommand{\abstractname}{\vspace{-\baselineskip}}

$$
\\\\
$$

\UDK{519.853.62}

\begin{center}
\textbf{УСКОРЕННЫЙ ГРАДИЕНТНЫЙ \aig{СЛАЙДИНГ-МЕТОД} В ЗАДАЧАХ МИНИМИЗАЦИИ СУММЫ ФУНКЦИЙ}

\textbf{Д.\,М.~Двинских,  C\,.C.~Омельченко, А.\,В.~Гасников, А.\,И.~Тюрин}


\end{center}

\begin{abstract}
\noindentВ статье \aig{предложен} новый способ обоснования ускоренного градиентного 
слайдинга \aig{Дж. Лана}, позволяющий распространить технику слайдинга на сочетание 
ускоренных градиентных методов с ускоренными методами редукции дисперсии. \aig{Получены новые оптимальные оценки} для решения задач минимизации суммы гладких сильно 
выпуклых функций с гладким регуляризатором.
\end{abstract}

\section{Введение} \label{section_1}


Многие задачи анализа данных (машинного обучения) \aig{приводят} к необходимости решения задач минимизации \aig{функции} вида суммы (эмпирический риск) с большим числом слагаемых, отвечающих объему выборки \cite{bib_1, bib_6, bib_18, bib_19}. В последнее десятилетие активно развиваются численные методы оптимизации функции вида суммы 
\cite{bib_1, bib_6, bib_8, bib_10}. В частности, были получены оптимальные методы (ускоренные методы редукции дисперсии) для такого класса задач, когда слагаемые в сумме гладкие (сильно) выпуклые функции, см., например, \cite{bib_10}. Были 
исследованы задачи, в которых дополнительно в \aig{минимизируемую функцию} вносится 
аддитивно, \aig{возможно, негладкий, но выпуклый / сильно выпуклый} композитный 
член (по терминологии анализа данных вносится слагаемое, отвечающее 
``регуляризации''), являющийся проксимально дружественным 
\cite{bib_11, bib_12}, т.е. задача минимизации такого члена с квадратичной добавкой -- простая задача. В 
настоящей работе предлагается способ получения оптимальных оценок для 
случая, когда композитный член будет выпуклым (сильно выпуклым) гладким, но 
уже не будет проксимально дружественным. Не предполагается проксимальная 
дружественность и у слагаемых в сумме.

В пункте \ref{section_2} техника ускоренного градиентного слайдинга Дж. Лана 
\cite[section 8.2]{bib_10} будет объяснена с помощью 
популярной в последнее время конструкции каталист 
\cite{bib_1, bib_13, bib_14}. Обнаруженный способ позволил распространить 
область приложений техники слайдинга на интересующий нас класс задач. В 
пункте 3 результаты пункта 2 обобщаются на различные негладкие постановки 
задач, в частности на обобщенные линейные модели 
\cite{bib_19} и другие модели, допускающие эффективное 
сглаживание \cite{bib_4, bib_17}.

\section{Основные результаты}\label{section_2}	
Рассмотрим следующую задачу
\begin{equation}
\label{eq1}
F\left( x \right)=f\left( x \right)+g\left( x \right)=f\left( x 
\right)+\frac{1}{m}\sum\limits_{k=1}^m {g_k \left( x \right)} \to \mathop 
{\min }\limits_x ,
\end{equation}
где $f$ и $g_k $ имеют $L_f $ и $L_g $-Липшицевы градиенты в 2-норме, а 
функция $F$ -- $\mu $-сильно выпуклая в 2-норме, причем $\mu \ll L_f $. \aig{В задаче \eqref{eq1} введем} дополнительное условие $m\le {L_g } \mathord{\left/ 
{\vphantom {{L_g } \mu }} \right. \kern-\nulldelimiterspace} \mu $. 
Результат Дж. Лана \cite[section 8.2]{bib_10} заключается в 
том, что для решения рассмотренной задачи с заданной точностью\footnote{ Не 
\aig{важно с} какой именно точностью $\varepsilon $. Эта точность будет входить 
под логарифмами в приведенные далее оценки, а для наглядности 
логарифмические сомножители было решено опустить. Далее оговорки о точности 
решения возникающих подзадач также опускаются, поскольку все это влияет 
только на логарифмические сомножители в итоговых оценках, которые опущены. 
Здесь и далее $\tilde {{\rm O}}\left( \right)={\rm O}\left( \right)$ с 
точностью до логарифмического множителя.} достаточно $\tilde {{\rm O}}\left( 
{\sqrt {{L_f } \mathord{\left/ {\vphantom {{L_f } \mu }} \right. 
\kern-\nulldelimiterspace} \mu } } \right)$ вычислений $\nabla f$ и $\tilde 
{{\rm O}}\left( {\sqrt {{L_g } \mathord{\left/ {\vphantom {{L_g } \mu }} 
\right. \kern-\nulldelimiterspace} \mu } } \right)$ вычислений $\nabla g$, 
т.е. $\tilde {{\rm O}}\left( {m\sqrt {{L_g } \mathord{\left/ {\vphantom 
{{L_g } \mu }} \right. \kern-\nulldelimiterspace} \mu } } \right)$ 
вычислений $\nabla g_k $.

Наложим еще одно дополнительное условие $\ag{mL_f \le L_g}$. Применим 
к рассмотренной задаче технику каталист \cite{bib_1, bib_13, bib_14}.\footnote{Заметим, 
что обойтись без этой техники не получается! Отметим также, что если 
использовать технику каталист в варианте \cite{bib_1, bib_9}, то применение данной техники не привносит 
дополнительного логарифмического множителя.} Тогда вместо исходной задачи 
(\ref{eq1}) потребуется $\tilde {{\rm O}}\left( {\sqrt {L \mathord{\left/ {\vphantom 
{L \mu }} \right. \kern-\nulldelimiterspace} \mu } } \right)$ раз решать 
задачу вида
\begin{equation}
\label{eq2}
f\left( x \right)+g\left( x \right)+\frac{L}{2}\left\| {x-x^k} \right\|_2^2 
\to \mathop {\min }\limits_x ,
\end{equation}
где $L$ по построению должно удовлетворять неравенству $\mu\le L \le L_f$. Задачу (\ref{eq2}) можно решать неускоренным композитным градиентным методом 
\cite{bib_1, bib_2, 
bib_3, bib_16}, считая $g\left( x 
\right)+\frac{L}{2}\left\| {x-x^k} \right\|_2^2 $ композитом. Число итераций 
такого метода будет совпадать с числом вычислений $\nabla f$ и равно $\tilde 
{{\rm O}}\left( {{L_f } \mathord{\left/ {\vphantom {{L_f } {\left( {L+\mu } 
\right)}}} \right. \kern-\nulldelimiterspace} {\left( {L+\mu } \right)}} 
\right)$. Но в условиях задачи не предполагалась проксимальная 
дружественность функции $g$, поэтому возникающую на каждой итерации 
неускоренного композитного градиентного метода задачу вида \aig{(детали см. в препринте \cite{bib_Ivanova})} 
\begin{equation}
\label{eq3}
\left\langle {\nabla f\left( {\tilde {x}^l} \right),x-\tilde {x}^l} 
\right\rangle +\frac{L_f }{2}\left\| {x-\tilde {x}^l} \right\|_2^2 +g\left( 
x \right)+\frac{L}{2}\left\| {x-x^k} \right\|_2^2 \to \mathop {\min 
}\limits_x ,
\end{equation}
в свою очередь, необходимо будет решать. Для решения задачи (\ref{eq3}) можно 
использовать ускоренный композитный метод редукции дисперсии 
\cite{bib_10, bib_11, bib_12}, считая $\frac{L_f}{2}\left\| {x-\tilde {x}^l} 
\right\|_2^2 +\frac{L}{2}\left\| {x-x^k} \right\|_2^2 $ композитом. Число 
вычислений $\nabla g_k $ для такого метода будет\footnote{\ag{Точнее говоря, оценка имеет вид:  $\tilde {\rm{O}}\left(m + \sqrt{mL_g/\left(L_f+L\right)}\right)$. Однако в виду предположений $mL_f \le L_g$, $L\le L_f$: $\tilde {\rm{O}}\left(m + \sqrt{mL_g/\left(L_f+L\right)}\right) = \tilde {\rm{O}}\left( \sqrt{mL_g/\left(L_f+L\right)}\right)$.}} $\tilde {{\rm O}}\left( 
{\sqrt {m{L_g } \mathord{\left/ {\vphantom {{L_g } {\left( {L_f +L+\mu } 
\right)}}} \right. \kern-\nulldelimiterspace} {\left( {L_f +L} 
\right)}} } \right)$. Таким образом, общее число вычислений $\nabla g_k $ 
будет\footnote{\ag{Первое слагаемое появилось из-за того, что в каталисте требуется считать $\nabla F$ на каждой итерации.} }
\begin{equation} \label{eq4}
   \ag{\tilde{\rm{O}}\left(m\sqrt{L/\mu}\right)+}\tilde{{\rm O}}\left( {\sqrt {L \mathord{\left/ {\vphantom {L \mu }} 
    \right. \kern-\nulldelimiterspace} \mu } } \right)\cdot \tilde {{\rm 
    O}}\left( {{L_f } \mathord{\left/ {\vphantom {{L_f } {\left( {L+\mu } 
    \right)}}} \right. \kern-\nulldelimiterspace} {\left( {L+\mu } \right)}} 
    \right)\cdot \tilde {{\rm O}}\left( {\sqrt {m{L_g } \mathord{\left/ 
    {\vphantom {{L_g } {\left( {L_f +L+\mu } \right)}}} \right. 
    \kern-\nulldelimiterspace} {\left( {L_f +L} \right)}} } \right).
\end{equation}
Выбирая параметр $L$ ($\mu\le L \le L_f$) так, чтобы выражение (\ref{eq4}) было минимальным, получим (с 
учетом сделанных предположений \ag{$mL_f \le L_g$} и $\mu \ll L_f )$, что 
$L\simeq L_f $. Следовательно, имеет место 

\begin{theorem}\label{th_1}

\textit{\ag{При $\ag{mL_f \le L_g}$} задачу (\ref{eq1}) можно решить с помощью описанной выше техники за }$\tilde {{\rm O}}\left( {\sqrt {{L_f } \mathord{\left/ 
{\vphantom {{L_f } \mu }} \right. \kern-\nulldelimiterspace} \mu } } 
\right)$\textit{ вычислений }$\nabla f$ и $\tilde {{\rm O}}\left( {\sqrt {m{L_g } \mathord{\left/ 
{\vphantom {{L_g } \mu }} \right. \kern-\nulldelimiterspace} \mu } } 
\right)$\textit{ вычислений }$\nabla g_k .$

\end{theorem}

Последняя оценка в $\tilde {{\rm O}}\left( {\sqrt m } \right)$ раз лучше оценки, которую можно получить, используя исходный ускоренный градиентный слайдинг Дж. Лана \cite{bib_10}. Несложно заметить \cite{bib_10}, что приведенные в теореме \ref{th_1} оценки оптимальны с точностью до логарифмических множителей.

Заметим, что в описанном выше подходе с \aig{$g(x)$ общего вида ускоренный метод редукции дисперсии можно заменить} на покоординантный спуск или безградиентный метод \cite{bib_7}. Таким образом, можно получить расщепление задачи не только по гладкости или структуре слагаемых,но и по структуре оракула, доступного для каждого из слагаемых. Другой пример такого расщепления см. в \cite{bib_5}.

Заметим также, что если в описанном выше подходе ограничиться вариантом каталиста из \cite{bib_13, bib_14}, то все рассуждения можно провести в модельной (для $f)$ общности \cite{bib_1, bib_3}.

		
\section{Приложение}

\aig{Заметим, что аналогично случаям задач из \cite{bib_15, bib_18, bib_19} описанная выше техника может использоваться и тогда, когда $g_k $ -- негладкие функции, но, 
допускающие, сглаживание \cite{bib_4, bib_17}. Скажем, двойственное сглаживание по Ю.Е. 
Нестерову \cite{bib_1, bib_4, bib_17}}. А именно, предположим, что функции $g_k $ 
имеют проксимально-дружественные сопряженные функции $g_k^\ast $. В 
частности, это имеет место для обобщенной линейной модели 
\cite{bib_19}, в которой $g_k \left( x \right):=g_k \left( 
{\left\langle {a_k ,x} \right\rangle } \right)$. Тогда, регуляризируя 
сопряженные функции $g_k^\ast $ с коэффициентом регуляризации $\sim 
\varepsilon $, где $\varepsilon $ -- желаемая точность (по функции) решения 
исходной задачи, \aig{получим, что} $\varepsilon 
\mathord{\left/ {\vphantom {\varepsilon 2}} \right. 
\kern-\nulldelimiterspace} 2$-решение сглаженной задачи будет $\varepsilon 
$-решением исходной. При том, что для сглаженной задачи $L_g \sim 
\varepsilon ^{-1}$.

Заметим, что с помощью регуляризации исходной задачи \cite{bib_1} описанные выше результаты распространяются с сильно выпуклого случая на просто выпуклый случай. Для этого в постановку выпуклой задачи \eqref{eq1} вносится регуляризация $+ \mu/2 \|x\|_2^2$, где $\mu = \varepsilon/R^2$. Здесь $\varepsilon$ -- желаемая точность решения \aig{задачи} по функции, а $R = \|x_*\|_2$ -- $2$-норма решения (на практике можно брать оценку сверху \aig{\cite{bib_1}}). Из \cite{bib_1} следует, что $\varepsilon/2$-решение так регуляризованной задачи будет $\varepsilon$ решением исходной задачи \eqref{eq1}. Продемонстрируем возможные преимущества предложенного подхода в выпуклом (но не сильно выпуклом случае).

Рассматривается постановка задачи
$$F(x) = \frac{1}{2}\langle x,Cx \rangle + \frac{1}{m}\sum_{k=1}^m g_k(\langle a_k, x \rangle) \to \min_{x\in\mathbb{R}^n}.$$
Предполагаем, что $|g_k''(y)| = O(1/\varepsilon)$, матрица $A = [a_1, ..., a_m]^T$ имеет $ms$ ненулевых элементов, $\max_{k=1,...,m} \|a_k\|_2^2 = O(s)$, где $1\ll s \le n$ и $C$ -- неотрицательно определенная матрица с $\lambda_{\max}(C)\le 1/(\varepsilon m)$. 
Ускоренный градиентный метод (FGM) \cite{bib_2} будет требовать $$O\left(\sqrt{\frac{\left(s/\varepsilon + \lambda_{\max}(C)\right)R^2}{\varepsilon}}\right)$$
итераций для достижения точности $\varepsilon$ по функции \aig{со сложностью} одной итерации
$$O\left(ms + n^2\right)$$ арифметических операций (а.о.).
\aig{В настоящей работе} предложен подход, который требует
$$\tilde{O}\left(\sqrt{\frac{ \lambda_{\max}(C)R^2}{\varepsilon}}\right)$$
итераций ускоренного градиентного метода для квадратичной формы (первого слагаемого). При этом сложность одной такой итерации
$$O(n^2)\quad \text{а.о.}$$
Также предложенный подход требует
$$\tilde{O}\left(\sqrt{\frac{\left(ms/\varepsilon \right)R^2}{\varepsilon}}\right)$$
итераций ускоренного метода редукции дисперсии  \cite{bib_1,bib_10,bib_11}. При этом сложность одной такой итерации
$$O(s) \quad \text{а.о.}$$

Для наглядности эти результаты собраны \aig{в таблицу 1}. \aig{Из таблицы 1} можно сделать вывод, что при $s\gg 1$, $\lambda_{\max}(C)\le 1/(\varepsilon m)\ll s/\varepsilon$, предложенный в данной \aig{работе} подход имеет лучшую теоретическую сложность, чем ускоренный градиентный метод, который принято было считать наилучшим для данного класса задач.

\begin{center}
\begin{table}[H]
\label{table_0}
\begin{tabular}{|c|c|c|ll}
\cline{1-3}
Алгоритм  & Сложность&  Ссылка&   &  \\ \cline{1-3}
FGM & $ O\left(\frac{R}{\varepsilon}\sqrt{s}\left(ms + n^2\right)\right) $ & \cite{bib_2} &  &  \\ \cline{1-3}
Слайдинг & $ \tilde{O}\left(\frac{R}{\varepsilon}\sqrt{ms}\cdot s\right) + \tilde{O}\left(\sqrt{\frac{\lambda_{\max}(C)R^2}{\varepsilon}}\cdot n^2\right) $ & данная статья &  &  \\ 
 \cline{1-3}
\end{tabular}
\caption{\aig{Сравнение алгоритмов}}
\end{table}
\end{center}



Работа поддержана грантами РФФИ 18-31-20005 мол\_а\_вед в п. 2 и РФФИ 19-31-90062 Аспиранты в п. 3.

Двинских Дарина Михайловна, Weierstrass Institute for Applied Analysis and Stochastics, 10117, Germany, Berlin, Mohrenstraße 39,  Московский физико-технический институт, 141701, Московская область, г. Долгопрудный, Институтский переулок, д.9.; Институт Проблем Передачи Информации РАН, 127994, г. Москва, ГСП-4, Большой Каретный переулок, 19, стр. 1; Тел. +7 (977) 365 16 24, e-mail: darina.dvinskikh@wias-berlin.de

Омельченко Сергей Сергеевич, Московский Физико-Технический Институт, 141701, Московская область, г. Долгопрудный, Институтский переулок, д.9; e-mail: sergey.omelchenko@phystech.edu.

Гасников Александр Владимирович, Московский Физико-Технический Институт, 141701, Московская область, г. Долгопрудный, Институтский переулок, д.9.; Институт Проблем Передачи Информации РАН, 127994, г. Москва, ГСП-4, Большой Каретный переулок, 19, стр. 1; Тел. +7 (905) 780 69 74, e-mail: gasnikov@yandex.ru.


 Тюрин Александр Игоревич, Национальный исследовательский университет Высшая школа экономики, 101000, г. Москва, ул. Мясницкая, 20;
 e-mail: alexandertiurin@gmail.com.



\newpage

\begin{thebibliography}{99}


\bibitem{bib_1}
\textit{ А.В.~Гасников}
 Современные численные методы оптимизации. Метод универсального градиентного
спуска. // arXiv:1711.00394

\bibitem{bib_2}
\textit{  А.В.~Гасников, Ю.Е.~Нестеров}
 Универсальный метод для задач стохастической композитной
оптимизации //  ЖВМ и МФ.
2018. V.58. № 1. P. 51--68
 

\bibitem{bib_3}
\textit{  А.В.~Гасников, А.И.~Тюрин}
 Быстрый градиентный спуск для задач выпуклой минимизации с оракулом, выдающим $\left( {\delta ,L} \right)$-модель функции в запрошенной точке//
 ЖВМ и МФ.
 2019. V.59. № 7. P. 1137--1150


\bibitem{bib_4}
\textit{  Z.~Allen-Zhu, E.~Hazan}
 Optimal black-box reductions between optimization objectives //
arXiv:1603.05642

\bibitem{bib_5}
\textit{  A.~Beznosikov, E.~Gorbunov, A.~Gasnikov}
 Derivative-free method for decentralized distributed non-smooth optimization//
arXiv:1911.10645

\bibitem{bib_6}
\textit{ L.~Bottou, F.E.~Curtis, J.~Nocedal}
 Optimization methods for large-scale machine learning//
arXiv:1606.04838

\bibitem{bib_7}
\textit{  P.~Dvurechensky, A.~Gasnikov, A.~Tiurin }
 Randomized Similar Triangles Method: A unifying framework for accelerated randomized optimization methods (Coordinate Descent, Directional Search, Derivative-Free Method)//
 arXiv:1707.08486

\bibitem{bib_8}
\textit{ E.~Hazan}
 Lecture notes: Optimization for Machine Learning//
arXiv:1909.03550

\aig{\bibitem{bib_Ivanova}
\textit{A.~Ivanova, A.~Gasnikov, P.~Dvurechensky, D.~Dvinskikh, A.~Tyurin, E.~Vorontsova, D.~Pasechnyuk}
Oracle Complexity Separation in Convex Optimization // arXiv:2002.02706}

\bibitem{bib_9}
\textit{  A.~Ivanova, D.~Grishchenko, A.~Gasnikov, E.~Shulgin}
 Adaptive Catalyst for smooth convex optimization //
arXiv:1911.11271

\bibitem{bib_10}
\textit{ G.~Lan}
 Lectures on optimization. Methods for Machine Learning.//
\url{https://pwp.gatech.edu/guanghui-lan/publications/}

\bibitem{bib_11}
\textit{  G.~Lan, Z.~Li, Y.~Zhou}
 A unified variance-reduced accelerated gradient method for convex optimization//
arXiv:1905.12412

\bibitem{bib_12}
\textit{  G.~Lan, Y.~Zhou}
 Randomized gradient extrapolation for distributed and stochastic optimization//
SIAM Journal on Optimization
2018. V.28. № 4. P. 2753--2782

\bibitem{bib_13}
\textit{  H.~Lin ,  J.~Mairal, Z.~Harchaoui }
 A universal catalyst for first-order optimization//
Proceedings of $29^{th}$ International conference Neural Information Processing Systems (NIPS).
2015.


\bibitem{bib_14}
\textit{  H.~Lin, J.~Mairal, Z.~Harchaoui}
 Catalyst acceleration for first-order convex optimization: from theory to practice//
arXiv:1712.05654


\bibitem{bib_15}
\textit{  K.~Mishchenko, P.~Richtarik}
 A Stochastic decoupling method for minimizing the sum of smooth and non smooth functions//
arXiv:1905.11535

\bibitem{bib_16}
\textit{ Yu.~Nesterov}
 Gradient methods for minimizing composite functions//
Math. Prog.
2013. V.140. № 1. P. 125--161 

\bibitem{bib_17}
\textit{ Yu.~Nesterov}
 Smooth minimization of non-smooth function//
Math. Program.
2005. V.103. № 1. P. 127--152 

\bibitem{bib_18}
\textit{ S.~Shalev-Shwartz, S.~Ben-David}
 Understanding Machine Learning: From theory to algorithms//
Cambridge University Press.
2014.

\bibitem{bib_19}
\textit{ S.~Shalev-Shwartz, O.~Shamir, N.~Srebro, K.~Sridharan} 
 Stochastic Convex Optimization//
COLT.
2009.

\bibitem{gasnikov016accrand}
\textit{A. Gasnikov, P. Dvurechensky, I. Usmanova} 
About accelerated randomized methods //
Proceedings of the Moscow Institute of Physics and Technology
2016.
V.8. №2 (30).


\bibitem{nesterov2012efficiency}
\textit{Y. Nesterov} 
Efficiency of coordinate descent methods on huge-scale optimization problems // SIAM Journal on Optimization 2012. V.22. № 2. P. 341--362

\bibitem{nesterov2017efficiency}
\textit{Y. Nesterov, S.U. Stich} 
\textit{Efficiency of the Accelerated Coordinate Descent Method on Structured Optimization Problems  //
SIAM Journal on Optimization
2017. V.27. № 1. P. 110--123}

\bibitem{fercoq2015accelerated}
\textit{O. Fercoq, P. Richt{\'a}rik} 
Accelerated, parallel, and proximal coordinate descent // SIAM Journal on Optimization 2015. V.25. № 4. P. 1997--2023
 

\end{thebibliography}

\end{document}